\documentclass[10pt]{amsart}
\usepackage{amsfonts,amstext,amssymb}
\newtheorem{theorem}{Theorem}
\newtheorem{claim}{Claim}
\newtheorem{lemma}{Lemma}
\newtheorem{remark}{Remark}
\newcommand{\powerinfront}[2]{{}^#1 #2}

\newcommand{\conn}{\operatorname{conn}}
\newcommand{\ran}{\operatorname{ran}}
\newcommand{\interior}{\operatorname{int}}
\newcommand{\id}{\operatorname{id}}

\newcommand{\0}{{\bf 0}}
\newcommand{\1}{{\bf 1}}
\newcommand{\lsup}{\sqcup}
\newcommand{\linf}{\sqcap}
\renewcommand{\land}{\wedge}
\renewcommand{\lor}{\vee}
\renewcommand\newsymbol[5]{%
\newcommand{\bigsqcap}{\setlength{\unitlength}{10pt}
                        \begin{picture}(1,1.25)
                            \thicklines
                            \put(0,-0.4){\line(0,1){1.4}}
                            \put(0,1){\line(1,0){1}}
                            \put(1,1){\line(0,-1){1.4}}
                        \end{picture}\;}
\DeclareMathSymbol#1{#3}%
   {\ifcase #2\or AMSa\or AMSb\fi}{"#4#5}}
\newsymbol\beperk 1216
\begin{document}
\begin{abstract}
    In this paper we give new proofs of the theorem of Ma\'{c}kowiak
    and Tymchatyn that every metric continuum is a
    weakly-confluent image of some one-dimensional hereditarily
    indecomposable continuum of countable weight. The first is a
    model-theoretic argument; the second is a topological proof
    inspired by the first.
\end{abstract}
\title{On the Ma\'{c}kowiak-Tymchatyn theorem}
\author{K. P. Hart}
\author{B. J. van der Steeg}
\address{Faculty of Information Technology and Systems\\
         TU Delft\\
         Postbus 5031\\
         2600~GA{}Delft\\
         the Netherlands}

\email[K. P. Hart]{k.p.hart@its.tudelft.nl}

\email[B.J. van der Steeg]{b.j.vandersteeg@its.tudelft.nl}

\urladdr[K. P.Hart]{http://aw.twi.tudelft.nl/\~{}hart}

\maketitle
%%%%%%%%%%%%%%%%%%%%%%%%%%%%%%%%%%%%%%%%%%%%%%%%%%%%%%%%%%%%%
%                                                           %
%                       INTRODUCTION                        %
%                                                           %
%%%%%%%%%%%%%%%%%%%%%%%%%%%%%%%%%%%%%%%%%%%%%%%%%%%%%%%%%%%%%
%
\section{Introduction}
In~\cite{MT} Ma\'{c}kowiak and Tymchatyn proved that every metric
continuum is the continuous image of a one-dimensional
hereditarily indecomposable continuum by a weakly confluent map.
In~\cite{HvMP} this result was extended to general continua, with
two proofs, one topological and one model-theoretic. Both proofs
made essential use of the metric result.

The original purpose of this paper was to (re)prove the metric
case by model-theoretic means. After we found this proof we
realized that it could be combined with any standard proof of the
completeness theorem of first-order logic (see e.g.,
Hodges~\cite{H}, 6.1]) to produce an inverse-limit proof of the
general form of the Ma\'{c}kowiak-Tymchatyn result. We present
both proofs. The model-theoretic argument occupies
sections~\ref{modprelim} and~\ref{modproof}, and the inverse-limit
approach appears in section~\ref{topproof}.

We want to take this opportunity to point out some connections
with work of Bankston~\cite{B}, who dualized the model-theoretic
notions of existentially closed structures and existential maps to
that of co-existentially closed compacta and co-existential maps.
He proves that co-existential maps are weakly confluent, that
co-existentially closed continua are one-dimensional and
hereditarily indecomposable, and that every continuum is the
continuous image of a co-existentially closed one. The map can in
general not be chosen co-existential, because co-existential maps
preserve indecomposability and do not raise dimension.
%
%%%%%%%%%%%%%%%%%%%%%%%%%%%%%%%%%%%%%%%%%%%%%%%%%%%%%%%%%%%%%
%                                                           %
%                       PRELIMINARIES                       %
%                                                           %
%%%%%%%%%%%%%%%%%%%%%%%%%%%%%%%%%%%%%%%%%%%%%%%%%%%%%%%%%%%%%
\section{Preliminaries}
\subsection{Ma\'{c}kowiak-Tymchatyn theorem}
The theorem of Ma\'{c}kowiak and~Tymchatyn, we are dealing with in
this paper states that every metric continuum is a weakly
confluent image of a one-dimensional hereditarily indecomposable
continuum of countable weight.

A continuum is \emph{decomposable} if it can be written as a union
of two proper subcontinua, it is called \emph{indecomposable} if
this is not the case. We call a continuum \emph{hereditarily
indecomposable} if every subcontinuum is indecomposable. This is
equivalent to saying that every two subcontinua that meet, one is
contained in the other. As in~\cite{HvMP} we can extend this
notion for arbitrary compact Hausdorff spaces. So a compact
Hausdorff space is hereditarily indecomposable if for every two
subcontinua that meet, one is contained in the other. We call a
continuous mapping between two continua \emph{weakly confluent} if
every subcontinuum in the range is the image of a subcontinuum in
the domain.
\begin{theorem}[Ma\'{c}kowiak and Tymchatyn \cite{MT}]
    \label{MaTym}
    Every metric continuum is a weakly confluent
    image of some one-dimensional hereditarily
    indecomposable continuum of the same weight.
\end{theorem}
In~\cite{HvMP} Hart, van Mill and Pol showed that the
Ma\'{c}kowiak and Tymchatyn result above implies the theorem for
the non-metric case using model-theoretic means.
\subsection{Wallman space}
In the proof we will consider the lattice of closed sets of our
metric continuum $X$ and try to find, through model-theoretic
means, another lattice in which we can embed our lattice of closed
sets of $X$. This new lattice will be a model for some sentences
which will make sure that its Wallman representation is a
continuum with certain properties. So at the base of the proof is
Wallman's generalization, to the class of distributive lattices,
of Stone's representation theorem for Boolean algebras. Wallmann's
representation theorem is as follows.
\begin{theorem}[\cite{W}]
\label{wallman} If $L$ is a distributive lattice, then there is a
compact $T_1$ space $X$ with a base for its closed sets that is a
homomorphic image of $L$. If $L$ is also disjunctive then we can
find a base for its closed sets that is an isomorphic image of
$L$. \end{theorem}
We call the space $X$ a Wallman space of $L$ or a Wallman
representation of $L$, notation: $wL$.

A lattice $L$ is \emph{disjunctive} if it models the sentence
\begin{equation}
\label{disjunctive} \forall a\,b\exists x[(a\linf
b\not=a)\rightarrow ((a\linf x=x)\land (b\linf x=\0))].
\end{equation}
Furthermore the space $X$ in theorem~\ref{wallman} is Hausdorff if
and only if the lattice $L$ is a normal lattice. We call a lattice
normal if it models the sentence
\begin{equation}
\label{normal} \forall a\,b\exists x\,y[(a\linf
b=\0)\rightarrow((a\linf x=\0)\land (b\linf y=\0)\land(x\lsup
y=\1))].
\end{equation}
Note that, if we start out with a compact Hausdorff space $X$ and
look at a base for its closed subsets which is closed under finite
unions and intersections, i.e., a (normal, disjunctive and
distributive) lattice, then the Wallman space of this lattice is
just the space $X$.
\begin{remark}
From now on we refer to a base for the closed subsets of some
topological space which is closed under finite unions and
intersections as a lattice base for the closed sets of the space
$X$.
\end{remark}
The following theorem shows how to create an onto mapping from
maps between lattices. In this theorem $2^X$ denotes the family of
all closed subsets of the space $X$.
\begin{theorem}\cite{DH}
\label{contimage}
    Let $X$ and $Y$ be compact Hausdorff spaces and let
    $\mathcal{C}$ be a base for the closed subsets of $Y$ that is
    closed under finite unions and intersections. Then $Y$ is a
    continuous image of $X$ if and only if there is a map
    $\phi:\mathcal{C}\rightarrow 2^X$ such that
    \begin{enumerate}
    \item
        $\phi(\emptyset)=\emptyset$, and if $F\not=\emptyset$ then
        $\phi(F)\not=\emptyset$
    \item
        if $F\cup G=Y$ then $\phi(F)\cup\phi(G)=X$
    \item
        if $F_1\cap\cdots\cap F_n=\emptyset$ then
        $\phi(F_1)\cap\cdots\cap\phi(F_n)=\emptyset$.
    \end{enumerate}
\end{theorem}
So $Y$ is certainly a continuous image of $X$ if there is an
embedding of some lattice base of the closed sets of $Y$ into
$2^X$.
%--------------------------------------------------------------
%
%   TRANSLATION OF PROPERTIES
%
%----------------------------------------------------------------
\subsection{Translation of properties}
Our model-theoretic proof of theorem~\ref{MaTym} will be as
follows. Given a metric continuum $X$, we will construct a lattice
$L$ such that some lattice base of $X$ is embedded into $L$, the
Wallman representation $wL$ of $L$ is a one-dimensional
hereditarily indecomposable continuum and that for every
subcontinuum in $X$ there exists a subcontinuum of $wL$ that is
mapped onto it.

For this we need to translate things like being hereditarily
indecomposable, being of dimension less than or equal to one and
being connected in terms of closed sets only.

To translate hereditarily indecomposability we use the following
characterization, due to Krasinkiewicz and Minc.
\begin{theorem}[Krasinkiewicz and Minc]
    A compact Hausdorff space is hereditarily indecomposable if and
    only if it is crooked between every pair of disjoint closed
    nonempty subsets.
\end{theorem}
Which the  authors translated in~\cite{HvMP} into terms of closed
sets only as follows.
\begin{theorem}\cite{HvMP}\label{HerIndec}
    A compact Hausdorff space $X$ is hereditarily indecomposable
    if and only if whenever four closed sets $C$, $D$, $F$ and~$G$
    in $X$ are given such that $C\cap D=C\cap G=F\cap D=\emptyset$
    one can write $X$ as the union of three closed sets $X_0$,
    $X_1$ and~$X_2$ such that $C\subset X_0$, $D\subset X_2$,
    $X_0\cap X_1\cap G=\emptyset$, $X_0\cap X_2=\emptyset$ and
    $X_1\cap X_2\cap F=\emptyset$.
\end{theorem}
So a compact Hausdorff space is hereditary indecomposable if the
lattice $2^X$ models the sentence
\begin{eqnarray}
    \lefteqn{\forall a\,b\,c\,d\ \exists
    x\,y\,z[((a\linf b=\0)\land (a\linf c=\0)
    \land (b\linf d=\0))\rightarrow}\label{HI}\\
    & & \rightarrow ((a\linf(y\lsup z)=\0)\land
    (b\linf (x\lsup y)=\0)\land (x\linf y=\0)\land\nonumber\\
    & & \land (x\linf y\linf d=\0)\land (y\linf z\linf
    c=\0)\land (x\lsup y\lsup z=\1))]\nonumber.
\end{eqnarray}
A space $X$ is of dimension less than or equal to one if the
lattice $2^X$ models the sentence
\begin{eqnarray}
    \lefteqn{\forall a\,b\,c\ \exists x\,y\,z[(a\linf b\linf
    c=\0)\rightarrow}\label{dim=1}\\
    & \rightarrow &((a\linf x=a)\land (b\linf y=b)\land
    (c\linf z=c)\land\nonumber\\
    & & \land (x\linf y \linf z=\0)\land (x\lsup y\lsup
    z=\1))]\nonumber.
\end{eqnarray}
A space $X$ is connected if the lattice $2^X$ models the sentence
$\conn(\1)$, where $\conn(a)$ is shorthand for the formula
$\forall x\, y[((x\linf y=\0)\land (x\lsup
y=a))\rightarrow(x=a)\lor (x=\0))]$.

\begin{remark}
For the next two sections, section~\ref{modprelim} and
section~\ref{modproof} we fix some metric continuum~$X$ and we
will show there exists a hereditarily indecomposable
one-dimensional continuum $Y$ of weight~$\operatorname{w}(X)$ such
that $X$ is weakly confluent image of~$Y$.
\end{remark}
%%%%%%%%%%%%%%%%%%%%%%%%%%%%%%%%%%%%%%%%%%%%%%%%%%%%%%%%%%%%%%%%%%
%                                                                %
% A CONTINUOUS IMAGE OF An HI 1-DIM CONTINUUM OF THE SAME WEIGHT %
%                                                                %
%%%%%%%%%%%%%%%%%%%%%%%%%%%%%%%%%%%%%%%%%%%%%%%%%%%%%%%%%%%%%%%%%%
\section{A continuous image of an hereditarily indecomposable
one-dimensional continuum of the same weight}\label{modprelim}

Using theorem~\ref{wallman} and~\ref{contimage} of the previous
section we see that to get a hereditarily indecomposable
one-dimensional continuum of weight $w(X)$ that maps onto $X$ we
must find a countable distributive, disjunctive normal lattice $L$
such that it is a model of the sentences~\ref{HI}, \ref{dim=1}
and~$\conn(\1)$, and furthermore some lattice base for the closed
sets of $X$ is embedded into this lattice~$L$.

Fix a lattice base $\mathcal{B}$ for the closed sets of $X$.

For some countable set of constants $K$ we will construct a set of
sentences $\Sigma$ in the language $\{\linf,\lsup,\0,\1\}\cup K$.
We will make sure that $\Sigma$ is a consistent set of sentences
such that, if we have a model $\mathfrak{A}=(A,\mathcal{I})$ for
$\Sigma$ then

$$L(\mathfrak{A})=\mathcal{I}\beperk K$$

is the universe of some lattice model in the language
$\{\linf,\lsup,\0,\1\}$ which is normal, distributive and
disjunctive and models the sentences~\ref{HI}, \ref{dim=1}
and~$\conn(\1)$. To make sure that $\mathcal{B}$ is embedded into
$L(\mathfrak{A})$ we simply add the diagram of the lattice
$\mathcal{B}$ to the set $\Sigma$ and make sure that there are
constants in $K$ representing the elements of $\mathcal{B}$. The
interpretations of $\linf$, $\lsup$, $\0$ and~$\1$ are given by
there interpretations under $\mathcal{I}$ in the model
$\mathfrak{A}$.

Let $K$ be the following countable set of constants
\begin{equation}
    \label{K}
    K=\bigcup_{-1\leq n<\omega}K_n=\bigcup_{-1\leq n<\omega}\{k_{n,m}:m<\omega\}.
\end{equation}
We will define the sentences of $\Sigma$ in an $\omega$-recursion.
So $\Sigma$ will be the set $\bigcup_{n<\omega}\Sigma_n$.

For definiteness we define $K_{-1}=\mathcal{B}$ and
$\Sigma_0=\triangle_{\mathcal{B}}$, the diagram of $\mathcal{B}$.

%----------------------------------------------------------------
%
%   CONSTRUCTION OF $\SIGMA$ IN \{\linf,\LSUP,\0,\1\}\CUP K$
%
%----------------------------------------------------------------
\subsection{Construction of $\Sigma$ in $\{\linf,\lsup,\0,\1\}\cup
K$}

Suppose we already defined the sentences up to $\Sigma_{5n}$.
\begin{enumerate}
    \item
    \label{Sigmalattice}
        $\Sigma_{5n+1}$ will be a set of sentences that will make sure
        that the supremum and infimum of any pair of constants in
        $\bigcup_{m\leq 5n}K_m$ are defined, using the new
        constants from $K_{5n+1}$.

        We will also make sure that the set
        $\Sigma_{5n+1}$ makes sure that distributivity holds for
        any triple of elements from $\bigcup_{m\leq 5n}K_m$.

        And we will make sure that the family of sets of sentences
        $\{\Sigma_{5n+1}:n<\omega\}$ will prevent the existence of a
        counterexample for $\conn(\1)$ in $\bigcup_{m\leq 5n}K_m$.
    \item
    \label{Sigmadisjunctive}
        $\Sigma_{5n+2}$ will be a set of sentences
        that will make sure that for every $a,b\in\bigcup_{m\leq 5n}K_m$ there
        exists some $c\in\bigcup_{m\leq n}K_{5m+2}$ such that the
        formula that is sentence~\ref{disjunctive} without
        quantifiers will hold for these $a$, $b$ and $c$.
    \item
    \label{Sigmanormal}
        $\Sigma_{5n+3}$ will be a set of sentences
        that will make sure that for every $a,b\in\bigcup_{m\leq
        5n}K_m$ there exist $c,d\in\bigcup_{m\leq n}K_{5m+3}$ such
        that the formula that is sentence~\ref{normal} without
        quantifiers will hold for these $a$, $b$, $c$ and~$d$.
    \item
        $\Sigma_{5n+4}$ will be a set of
        sentences that will make sure that the according to the
        elements of $K_{5n+4}\cup\bigcup_{m\leq 5n}K_m$ the
        dimension of the Wallman space of $L(\mathfrak{A})$
        for any model $\mathfrak{A}$ of $\Sigma$ will be
        less than or equal to one.
    \item
    \label{SigmaHI}
        $\Sigma_{5(n+1)}$ will be a set of
        sentences that will make sure that for any
        $a,b,c\in\bigcup_{m<5n}K_m$ there exist
        $x,y,z\in K_{5(n+1)}$ such that the formula, which is the
        sentence~\ref{HI} without quantifiers, holds for this
        $a,b,c$ and $x,y,z$.
\end{enumerate}
We now show how to define the sets of sentences of
$\{\linf,\lsup,\0,\1\}\cup\bigcup_{m<5n+4}K_m$ as described
in~\ref{Sigmalattice}~-~\ref{SigmaHI}.

We have a natural order $\triangleleft$ on the set $K=\bigcup_m
K_m$ defined by $$k_{n,m}\triangleleft k_{r,t}\leftrightarrow
[(n<r)\lor ((n=r)\land (m<t))].$$

Let $\{p_l\}_{l<\omega}$ be an enumeration of $$\{p\in
[\bigcup_{m\leq 5n} K_m]^2: p\setminus\bigcup_{m\leq 5(n-1)}
K_m\not=\emptyset\}.$$
\begin{eqnarray*}
\Sigma_{5n+1}^0 &=& \{\bigsqcap p_l=k_{5n+1,2l}:l<\omega\}\\
\Sigma_{5n+1}^1 &=& \{\bigsqcup p_l=k_{5n+1,2l+1}:l<\omega\}\\
\Sigma_{5n+1}^2 &=& \{a\lsup a=a, a\linf a=a: a\in\bigcup_{m\leq
    5n}K_m\}\\
\Sigma_{5n+1}^3 &=& \{a\lsup(b\lsup c)=(a\lsup b)\lsup c,
    a\linf(b\linf c)=(a\linf b)\linf c: a,b,c\in\bigcup_{m\leq
    5n}K_m\}\\
\Sigma_{5n+1}^4 &=& \{a\lsup(b\linf c)=(a\lsup b)\linf (a\lsup c):
    a,b,c\in\bigcup_{m\leq 5n}K_m\}\\
\Sigma_{5n+1}^5 &=& \{a\lsup(a\linf b)=a, a\linf(a\lsup b)=a:
    a,b\in\bigcup_{m\leq 5n}K_m\}\\
\Sigma_{5n+1}^6 &=& \{[((a\lsup b=\1)\land (a\linf
    b=\0))\rightarrow ((a=\0)\lor (a=\1))]:
    a,b\in\bigcup_{m\leq 5n}K_m\}
\end{eqnarray*}

Define $\Sigma_{5n+1}$ by $$\Sigma_{5n+1}:=\bigcup_{i\leq
6}\Sigma_{5n+1}^i.$$ This set of sentences will make sure that any
model of $\Sigma$ in the language $\{\linf,\lsup,\0,\1\}\cup K$
will be a distributive lattice and also a model of the sentence
$\conn(\1)$.
\begin{eqnarray*}
    \Sigma_{5n+2} &=& \{[(\max_\triangleleft
    p_l\linf\min_\triangleleft
    p_l=\0)\rightarrow
    ((\max_\triangleleft p_l\linf k_{5n+2,2l}=\0)\land\\
    & & \land (\min_\triangleleft p_l\linf k_{5n+2,2l+1}=\0)\land
    (k_{5n+2,2l}\lsup
    k_{5n+2,2l+1}=\1))]: l<\omega\}
\end{eqnarray*}
This set of sentences will make sure that any (lattice) model of
$\Sigma$ in the language $\{\linf,\lsup,\0,\1\}\cup K$ will be
normal.

The following set of sentences makes sure that any model of
$\Sigma$ in the language $\{\linf,\lsup,\0,\1\}\cup K$ which is
also a lattice is a disjunctive lattice.
\begin{eqnarray*}
    \Sigma_{5n+3}^0 = \{[(\max_\triangleleft p_l\linf
    \min_\triangleleft p_l\not=\max_\triangleleft p_l)
    & \rightarrow & ((k_{5n+3,2l}\linf\max_\triangleleft
    p_l=k_{5n+3,2l})\land\\
    & & \land (k_{5n+3,2l}\linf\min_\triangleleft
    p_l=\0))]:l<\omega\}\\
    \Sigma_{5n+3}^1=\{[(\min_\triangleleft
    p_l\linf\max_\triangleleft p_l\not=\min_\triangleleft p_l)
    &\rightarrow & ((k_{5n+3,2l+1}\linf\min_\triangleleft
    p_l=k_{5n+3,2l+1})\land\\
    & & \land (k_{5n+3,2l+1}\linf\max_\triangleleft
    p_l=\0))]:l<\omega\}
\end{eqnarray*}
And define $\Sigma_{5n+3}$ by
$$\Sigma_{5n+3}=\Sigma_{5n+3}^0\cup\Sigma_{5n+3}^1.$$

Let $\zeta$ denote the following formula in
$\{\linf,\lsup,\0,\1\}$
\begin{eqnarray*}
    \zeta(a,b,c;x,y,z) = [(a\linf b\linf c=\0) &\rightarrow &
    ((a\linf x=a)\land (b\linf y=b)\land (c\linf z=c)\land\\
    & & \land (x\linf y\linf z=\0)\land (x\lsup y\lsup z=\1))]
\end{eqnarray*}
Let $\{q_l\}_{l<\omega}$ be an enumeration of the set
\begin{equation*}
    \{q\in [\bigcup_{m\leq 5n}K_m]^3:q\setminus
    \bigcup_{m\leq 5(n-1)}K_m\not=\emptyset\}
\end{equation*}
For every $l<\omega$ write $q_l=\{q_l(0),q_l(1),q_l(2)\}$.

Now define $\Sigma_{5n+4}$ by
\begin{equation*}
    \Sigma_{5n+4}=\{\zeta(q_l(0),q_l(1),q_l(2);k_{5n+4,3l},
    k_{5n+4,3l+1},k_{5n+4,3l+2}):l<\omega\}.
\end{equation*}
This will make sure that the Wallman space of any lattice model of
$\Sigma$ will be at most one-dimensional.

For making sure that the Wallman space of any model of $\Sigma$
will be hereditarily indecomposable we introduce the following
formulas in the language $\{\linf,\lsup,\0,\1\}$:
\begin{eqnarray}
    \phi(a,b,c,d) & = & [(a\linf b=\0)\land (a\linf d=\0)\land
    (b\linf c=\0)]\nonumber\\
    \psi(a,b,c,d;x,y,z) & = & [(x\lsup y\lsup z=\1)\land
    (x\linf z=\0)\land\nonumber\\
    & & \land (a\linf (y\lsup z)=\0)\land (b\linf(x\lsup y)=\0)
    \land\nonumber\\
    & & \land (x\linf y\linf d=\0)\land (y\linf z\linf
    c=\0)]\nonumber\\
    \theta(a,b,c,d;x,y,z)  & = & \phi(a,b,c,d) \rightarrow
    \psi(a,b,c,d;x,y,z) \label{theta}
\end{eqnarray}
Let $\{r_l\}_{l<\omega}$ be an enumeration of the set $$\{r\in
\powerinfront{4}{[\bigcup_{m\leq 5n}K_m]}:
\ran(r)\setminus\bigcup_{m\leq 5(n-1)}K_m\not=\emptyset\}.$$

Let $\Sigma_{5(n+1)}$ be the set of sentences defined by:
\begin{equation*}
    \Sigma_{5(n+1)}=\{\theta(r_l(0),r_l(1),r_l(2),r_l(3);
    k_{5(n+1),3l}, k_{5(n+1),3l+1},k_{5(n+1),3l+2}):l<\omega\}
\end{equation*}
Here the formula $\theta$ is as in equation~\ref{theta}.
%----------------------------------------------------------------
%
%   CONSISTENCY OF $\SIGMA$ IN $\{\linf,\LSUP,\0,\1\}\CUP K$
%
%----------------------------------------------------------------
\subsection{Consistency of $\Sigma$ in $\{\linf,\lsup,\0,\1\}\cup
K$}\label{CON(Sigma)}

In this section we show that $\Sigma$ is a consistent set of
sentences.

We will find for $\Sigma'\in[\Sigma]^{<\omega}$ a metric space
$X(\Sigma')$ and an interpretation function
$\mathcal{I}:K\rightarrow 2^{X(\Sigma')}$ such that
$(X(\Sigma'),\mathcal{I})\models
\Sigma'\cup\triangle_\mathcal{B}$. The interpretations of $\linf$,
$\lsup$, $\0$ and~$\1$ will always be $\cap$, $\cup$ (the normal
set intersection and union), $\emptyset$ and~$X(\Sigma')$
respectively.

For $\Sigma'=\emptyset$ we let $X(\emptyset)=X$ and we interpret
every constant from $K_{-1}$ as its corresponding base element in
$\mathcal{B}$. Extend the interpretation function by assigning the
empty set to all constants of $K\setminus K_{-1}$. It is obvious
that
$(2^{X(\emptyset)},\mathcal{I})\models\triangle_{\mathcal{B}}$.
\begin{remark}
\label{wellorder} As the interpretation of $\linf$ and $\lsup$ in
the metric continuum $X(\Sigma')$ will always be the normal set
intersection and set union, all the sentences in $\Sigma_{5n+1}^i$
for some $n<\omega$ and $i\in\{3,4,5,6\}$ are true in the model
$(2^{X(\Sigma')},\mathcal{I})$. So we can ignore these sentences
and for the remainder of this section concentrate on the remaining
sentences of $\Sigma$.
\end{remark}
We can define a well order $\sqsubset$ on the set
$\Sigma\setminus\{\Sigma_{5n+1}^i:n<\omega\ \text{and}\
i\in\{3,4,5,6\}\}$ by stating that $\phi\sqsubset\psi$ if and only
if there are $n<m<\omega$ such that $\phi\in\Sigma_n$ and
$\psi\in\Sigma_m$ or there are $k<l<\omega$ and $n<\omega$ such
that $\phi,\psi\in\Sigma_n$ and $\phi$ is a sentence that mentions
$p_k$ ($q_k$ or $r_k$ respectively) and $\psi$ is a sentence that
mentions $p_l$ ($q_l$ or $r_l$ respectively).

Suppose $\Sigma'$ is a finite subset of $\Sigma$ such that for all
of its proper subsets $\Sigma''$  there exists a metric continuum
$X(\Sigma'')$ and an interpretation function
$\mathcal{I}:K\rightarrow 2^{X(\Sigma'')}$ such that
$(X(\Sigma''),\mathcal{I})\models
\Sigma''\cup\triangle_\mathcal{B}$.

Let $\theta$ be the $\sqsubset$-maximal sentence in
$\Sigma'\setminus\{\Sigma_{5n+1}^i:n<\omega\ \text{and}\
i\in\{3,4,5,6\}\}$. We will show that there exists a metric space
$X(\Sigma')$ and an interpretation function
$\mathcal{I}:K\rightarrow 2^{X(\Sigma')}$ such that
$(X(\Sigma'),\mathcal{I})\models\Sigma'\cup\triangle_\mathcal{B}$.

Let $\Sigma''=\Sigma'\setminus\{\theta\}$.
%-----------------------------------------------------------------
%-----------------------------------------------------------------
\subsubsection{$\theta\in\bigcup_{m<\omega}\{\Sigma_{5n+1}\cup
                \Sigma_{5m+2}\cup\Sigma_{5m+3}\}$}
\label{theta123} We can simply let $X(\Sigma')=X(\Sigma'')$ and
either (re)interpret the new constant as the intersection or union
of two closed sets in $X(\Sigma'')$ if $\theta$ is in some
$\Sigma_{5n+1}$ or, if $\theta$ is an element of some
$\Sigma_{5m+2}$ or $\Sigma_{5m+3}$, using the fact that the space
$X(\Sigma'')$ is normal find (re)interpretations for the newly
added constants, in an obvious way.
%-----------------------------------------------------------------
%-----------------------------------------------------------------
\subsubsection{$\theta\in\{\Sigma_{5m+4}:m<\omega\}$}
\label{theta4} Suppose the preamble of $\theta$ is true in the
model $(2^{X(\Sigma'')},\mathcal{I})$, where $\theta$ is the
following sentence
\begin{eqnarray*}
    \theta &=& [(a\linf b\linf c=\0)\rightarrow (a\linf x=
    a)\land (b\linf y=b)\land \\
    & &
    \land (c\linf z=c)\land (x\linf y\linf z=\0)\land (x\lsup y\lsup
    z=\1)].
\end{eqnarray*}
If $a$ has a zero interpretation then we can choose $x=\0$, $y=\1$
and $z=\1$, and this interpretation of $x$, $y$ and~$z$ makes sure
that $\theta$ holds in the model $(2^{X(\Sigma'')},\mathcal{I})$.
So we may assume that $a$, $b$ and~$c$ have non zero
interpretations.

As the space $X(\Sigma'')$ is metric, we can assume that we have a
metric $\rho$ on $X(\Sigma'')$. Moreover we can assume that $\rho$
is bounded by~$1$.

Consider the following function $f$ from $X(\Sigma'')$ to
$\mathbb{R}^3$

$$f(x)=(\kappa_a(x),\kappa_b(x),\kappa_c(x)),$$

where $\kappa_a:X(\Sigma'')\rightarrow [0,1]$ is defined by

$$\kappa_a(x)=\frac{\rho(x,a)}{\rho(x,a)+\rho(x,b)+\rho(x,c)},$$

and $\kappa_b$ and $\kappa_c$ are like $\kappa_a$, but with $a$
interchanged with $b$ and $c$ respectively. Then $f[X(\Sigma'')]$
is a subset of the triangle
$T=\{(t_1,t_2,t_3)\in\mathbb{R}^3:t_1+t_2+t_3=1\ \text{and}\
t_1,t_2,t_3\geq 0\}$.

The space $X(\Sigma'')$ is embedded in the space
$X(\Sigma'')\times T$ by the graph of $f$ (in other words the
embedding is defined by $x\mapsto (x,f(x))$). Let us denote this
embedding by $g$.

Consider the space $\partial T\times [0,1]$, where $\partial
T=T\setminus\interior(T)$ in $\mathbb{R}^3$. Let $h$ be the map
from $\partial T\times [0,1]$ onto $T$ defined by
$$h((x,t))=x(1-t)+t(\frac{1}{3},\frac{1}{3},\frac{1}{3}).$$ The
map $h$ restricted to $\partial T\times [0,1)$ is a homeomorphism
between $\partial T\times [0,1)$ and
$T\setminus\{(\frac{1}{3},\frac{1}{3},\frac{1}{3})\}$.

We define $X(\Sigma')$ as the space $$X(\Sigma')=(\id\times
h)^{-1}[g[X(\Sigma'')]].$$

Let us (re)interpret the constants $k$ in $K$ in the following
way: $$\mathcal{I}(k):=\mathcal{I}(k)\times(\partial
T\times[0,1])\cap X(\Sigma')\ (=(\id\times
h)^{-1}[g[\mathcal{I}(k)]]).$$

\begin{remark}
\label{mono-cl} For future reference we note that, as the inverse
images of points $(x,(t_1,t_2,t_3))$ under the map $\id\times h$
are points for $(x,(t_1,t_2,t_3))$ in $X(\Sigma'')\times T$ for
which $(t_1,t_2,t_3)\not=(\frac{1}{3},\frac{1}{3},\frac{1}{3})$
and equal to $\{x\}\times\partial T\times\{1\}$ for those
$(x,(t_1,t_2,t_3))$ in $X(\Sigma'')\times T$ for which
$(t_1,t_2,t_3)=(\frac{1}{3},\frac{1}{3},\frac{1}{3})$, we have
that the map $\id\times h:X(\Sigma'')\times(\partial
T\times[0,1])\rightarrow X(\Sigma'')\times T$ is monotone.
Furthermore it is also closed.
\end{remark}

We did nothing to disturb the truth or falsity of the sentences
$\Sigma''$ in the model $(2^{X(\Sigma'')},\mathcal{I})$ as
$f^{-1}[A\cap B]=f^{-1}[A]\cap f^{-1}[B]$ and $f^{-1}[A\cup
B]=f^{-1}[A]\cup f^{-1}[B]$ for any function $f$ and any sets $A$
and $B$.

So we have that $(2^{X(\Sigma')},\mathcal{I})$ is a model for
$\Sigma''$.

Let $A$ be the line segments between $(0,1,0)$ and $(0,0,1)$, $B$
the line segment between $(1,0,0)$ and $(0,0,1)$ and $C$ the line
segment between $(1,0,0)$ and $(0,1,0)$. Now we (re)interpret $x$,
$y$ and~$z$ as follows
\begin{eqnarray*}
\mathcal{I}(x) &:=& X(\Sigma'')\times(A\times [0,1])\cap
X(\Sigma')\\ \mathcal{I}(y) &:=& X(\Sigma'')\times(B\times
[0,1])\cap X(\Sigma')\\ \mathcal{I}(z) &:=&
X(\Sigma'')\times(C\times [0,1])\cap X(\Sigma')
\end{eqnarray*}
As is easily seen, this interpretation of the constants $x$, $y$
and~$z$ makes the sentence $\theta$ a true sentence in the model
$(2^{X(\Sigma')},\mathcal{I})$. So
$(2^{X(\Sigma')},\mathcal{I})\models\Sigma'$.
%-----------------------------------------------------------------
%-----------------------------------------------------------------
\subsubsection{$\theta\in\{\Sigma_{5(m+1)}:m<\omega\}$}
\label{theta5} Suppose the preamble of $\theta$ is true in the
model $(2^{X(\Sigma'')},\mathcal{I})$, where $\theta$ is the
sentence $$\theta=\phi(a,b,c,d)\rightarrow\psi(a,b,c,d;x,y,z),$$
as in equation~\ref{theta}.

If the interpretation of $a$ is zero we can simply take $x=y=\0$
and $z=\1$ to make $(2^{X(\Sigma'')},\mathcal{I})$ a model of
$\theta$. So again we may assume that the interpretations of $a$,
$b$, $c$ and~$d$ are nonzero.

To show that $\Sigma'$ is a consistent set of sentences we are
going to use an idea from~\cite{HvMP}.

With the aid of Urysohn's lemma we can find a continuous function
$f:X(\Sigma'')\rightarrow [0,1]$ such that
$f(\mathcal{I}(a))\subset\{0\}$, $f(\mathcal{I}(b))\subset\{1\}$,
$f(\mathcal{I}(c))\subset [0,\frac{1}{2}]$ and
$f(\mathcal{I}(d))\subset[\frac{1}{2},1]$.

Let $P$ denote the (closed and connected) subset of
$[0,1]\times[0,1]$ given by
$$P=\{\frac{1}{4}\}\times[0,\frac{2}{3}]\cup
[\frac{1}{4},\frac{1}{2}]\times\{\frac{2}{3}\}\cup
\{\frac{1}{2}\}\times[\frac{1}{3},\frac{2}{3}]\cup
[\frac{1}{2},\frac{3}{4}]\times\{\frac{1}{3}\}\cup
\{\frac{3}{4}\}\times[\frac{1}{3},1].$$
Let $X^+\subset [0,1]\times X(\Sigma'')$ denote the pre-image of
the set $P$ under the function $\id\times f$:

$$X^+=\{(t,x)\in [0,1]\times X(\Sigma''):(t,f(x))\in P\}.$$

As $P$ is closed and $\id\times f$ is continuous we have that
$X^+$ is a compact metric space. Define the (continuous) map
$\pi:X^+\rightarrow X(\Sigma'')$ by $\pi((t,x))=x$ for every
$(t,x)\in X^+$.
\begin{lemma}
    There exists a unique component $C$ of $X^+$ such that $\pi[C]=X(\Sigma'')$.
    \label{cpnt_in_X+}
\end{lemma}
\begin{proof}
    Suppose we have closed sets $F$ and $G$ such that $X^+=F+G$.
    Define subsets $A_i,B_i$ of $X$, where $i\in\{0,1,2\}$ , by
    \begin{eqnarray*}
        A_0 & = & \{x\in X(\Sigma''):(\frac{1}{4},x)\in F\},\
        B_0 = \{x\in X(\Sigma''):(\frac{1}{4},x)\in G\}\\
        A_1 & = & \{x\in X(\Sigma''):(\frac{1}{2},x)\in F\},\
        B_1 = \{x\in X(\Sigma''):(\frac{1}{2},x)\in G\}\\
        A_2 & = & \{x\in X(\Sigma''):(\frac{3}{4},x)\in F\},\
        B_2 = \{x\in X(\Sigma''):(\frac{3}{4},x)\in G\}
    \end{eqnarray*}
    It is clear that $A_i\cap B_i=\emptyset$ for every
    $i\in\{0,1,2\}$.
    \begin{claim}
        The following holds
        \begin{enumerate}
        \item
            For every $x\in(A_0\cap B_1)\cup (B_0\cap A_1)$ we have
            $f(x)<\frac{2}{3}$.
        \item
            For every $x\in(A_1\cap B_2)\cup (B_1\cap A_2)$ we have
            $f(x)>\frac{1}{3}$.
        \end{enumerate}
    \end{claim}
    \begin{proof}
        As the proofs of the statements are very similar we will
        only prove the first statement.

        If $x\in A_0\cap B_1$ or $x\in B_0\cap A_1$ then
        $f(x)\leq\frac{2}{3}$. As $f(x)=\frac{2}{3}$ is impossible,
        we are done.
    \end{proof}
    Let us define the following closed sets $A^*$ and~$B^*$ of
    $X(\Sigma'')$ by
    \begin{eqnarray*}
        A^* = &\bigcup&\{f^{-1}[0,\frac{1}{3}]\cap
        A_0,f^{-1}[\frac{2}{3},1]\cap A_2,A_0\cap A_1\cap A_2,\\
            & &A_0\cap B_1\cap B_2,B_0\cap B_1\cap A_2\}\\
        B^*= &\bigcup&\{f^{-1}[0,\frac{1}{3}]\cap
        B_0,f^{-1}[\frac{2}{3},1]\cap B_2,B_0\cap B_1\cap B_2,\\
            & &B_0\cap A_1\cap A_2,A_0\cap A_1\cap B_2\}
    \end{eqnarray*}
    The sets $A^*$ and $B^*$ are disjoint closed subsets of $X(\Sigma'')$
    and their union is the whole of $X(\Sigma'')$. As $X(\Sigma'')$ is connected
    one of these sets must be empty. So without loss of generality
    we can assume that $B^*=\emptyset$.

    We see now that $\pi[F]=X(\Sigma'')$ and that
    $\pi[G]\subset f^{-1}[\frac{1}{3},\frac{2}{3}]$. It follows
    that whenever $C$ is a clopen subset of $X^+$ then
    either $\pi[C]=X(\Sigma'')$ and $\pi[X^+\setminus
    C]\subset f^{-1}[\frac{1}{3},\frac{2}{3}]$ or it is the other way
    around. This shows that $\mathcal{F}=\{C:C\ \text{is clopen
    and}\ \pi[C]=X(\Sigma'')\}$ is an ultrafilter in the family of
    clopen subsets of $X^+$; its intersection
    $\bigcap\mathcal{F}$ is the unique component of $X^+$
    that is mapped onto $X(\Sigma'')$. We let $C$ be this
    component. This ends the proof of lemma~\ref{cpnt_in_X+}.
\end{proof}
Let $X(\Sigma')$ be the unique component $C$ of $X^+$ that is
mapped onto $X(\Sigma'')$ by the map $\pi$ with the subspace
topology.

In $2^{X(\Sigma')}$, the constants $x$, $y$ and~$z$ that will make
the sentence $\theta$ true will have the following
(re)interpretations:
\begin{eqnarray*}
    \mathcal{I}(x) &=& \{(t,x)\in X(\Sigma'):t\in[0,\frac{3}{8}]\},\\
    \mathcal{I}(y) &=& \{(t,x)\in X(\Sigma')
    :t\in[\frac{3}{8},\frac{5}{8}]\}\ \text{and}\\
    \mathcal{I}(z) &=& \{(t,x)\in X(\Sigma'):t\in[\frac{5}{8},1]\}.
\end{eqnarray*}

The (re)interpretation of the constants in $K$ will be as follows.
$$\mathcal{I}(k):=[0,1]\times\mathcal{I}(k)\cap X(\Sigma')\
(=\pi^{-1}[\mathcal{I}(k)]\cap X(\Sigma')).$$

As $\pi$ maps $C$ onto $X(\Sigma'')$ we have that
$(2^{X(\Sigma')},\mathcal{I})$ is a model of $\Sigma'$, as the
truth or falsity of sentences in $\Sigma''$ are not affected by
the new interpretation of the constants.
%%%%%%%%%%%%%%%%%%%%%%%%%%%%%%%%%%%%%%%%%%%%%%%%%%%%%%%%%%%%%
%                                                           %
%                   THE Ma\'{c}kowiak-TymchATYN THEOREM      %
%                                                           %
%%%%%%%%%%%%%%%%%%%%%%%%%%%%%%%%%%%%%%%%%%%%%%%%%%%%%%%%%%%%%
\section{The Ma\'{c}kowiak-Tymchatyn theorem}
\label{modproof}
Apart from the weakly confluent property of the
continuous map we have proven the Ma\'{c}kowiak-Tymchatyn theorem,
theorem~\ref{MaTym}. To make sure that the continuous map
following from the previous section is  weakly confluent, we must
consider all the subcontinua of the space $X$.

We let $\hat{K}$ be the following set

\begin{eqnarray*}
\hat{K}=\bigcup_{-2\leq n<\omega}\hat{K}_n= \bigcup_{-2\leq
n<\omega}\{k_{n,\alpha}:\alpha<\vert 2^X\vert\}.
\end{eqnarray*}

We will construct a set

$$\hat{\Sigma}=\bigcup_{-1\leq n<\omega}\hat{\Sigma}_n$$

of sentences in the language $\{\linf,\lsup,\0,\1\}\cup\hat{K}$
similar as in the previous section such that given any model
$\mathfrak{A}=(A,\mathcal{I})$ of $\hat{\Sigma}$, the set
$L(\mathfrak{A})=\mathcal{I}\beperk\hat{K}$ will be the universe
of some normal distributive and disjunctive lattice such that
\begin{enumerate}
\item
    $L(\mathfrak{A})$ is a model of the sentences~\ref{HI},
    \ref{dim=1} and~$\conn(\1)$,
\item
    the lattice $2^X$ is embedded into  $L(\mathfrak{A})$ so there
    exists a continuous map $f$ from $wL(\mathfrak{A})$
    onto $X$,
\item
    for every subcontinuum of $X$ there exists a subcontinuum of
    $wL(\mathfrak{A})$ that is mapped onto it by $f$.
\end{enumerate}

%----------------------------------------------------------------
%
%   A WEAKLY CONFLUENT MAP
%
%----------------------------------------------------------------
\subsection{A weakly confluent map}

We let $\hat{K}_{-1}=\{k_{-1,\alpha}<\vert 2^X\vert\}$ correspond
to the set $2^X=\{x_\alpha:\alpha<\vert 2^X\vert\}$ in such a way
that the set $\mathcal{C}(X)$ of all the subcontinua of $X$
corresponds to the set $\{x_\alpha:\alpha<\beta\}$ for some
ordinal number $\beta<\vert 2^X\vert$. Let the set of sentences
$\hat{\Sigma}_0$ in $\{\linf,\lsup,\0,\1\}\cup\hat{K}_{-1}$
correspond to $\triangle_{2^X}$, the diagram of the lattice $2^X$.

We want to define a set of sentences $\Sigma_{-1}$ in
$\{\linf,\lsup,\0,\1\}\cup K_{-2}\cup K_{-1}$ that will make sure
that if $\mathfrak{A}$ is a model of $\hat{\Sigma}$ in the
language $\{\linf,\lsup,\0,\1\}\cup\hat{K}$ then we have for every
subcontinuum in $X$ a subcontinuum of $wL(\mathfrak{A})$ that will
be mapped onto it by the continuous onto map we get by the fact
that $2^X$ is embedded in the lattice $L(\mathfrak{A})$.
\begin{eqnarray*}
    \hat{\Sigma}^0_{-1}&=&
    \{\conn(k_{-2,\alpha})\land (k_{-2,\alpha}\linf
    k_{-1,\alpha}=k_{-2,\alpha})):\alpha<\beta\}\\
    \hat{\Sigma}^1_{-1}&=&\{(\conn(k_{-2,\alpha})\land
    (k_{-2,\alpha}\linf k_{-1,\gamma}=k_{-2,\alpha}))\rightarrow\\
    & &\rightarrow (k_{-1,\alpha}\linf k_{-1,\gamma}=k_{-1,\alpha}):
    \alpha<\beta,\ \gamma<\vert 2^X\vert\}\\
    \hat{\Sigma}^2_{-1}&=&\{k_{-2,\gamma}=\0:\beta\leq\gamma<\vert
    2^X\vert\}.
\end{eqnarray*}
And define the set of sentences $\Sigma_{-1}$ as
\begin{equation}
    \label{weakly-confluent}
    \hat{\Sigma}_{-1}=\hat{\Sigma}^0_{-1}\cup\hat{\Sigma}^1_{-1}
    \cup\hat{\Sigma}^2_{-1}.
\end{equation}
Suppose $\mathfrak{A}$ is a model of $\hat{\Sigma}$. The set
$\hat{\Sigma}^0_{-1}$ will make sure that for every subcontinuum
$C$ of $X$ there is some subcontinuum $C'$ of $wL(\mathfrak{A})$
that is mapped into $C$ by the continuous onto map $f$ we get from
theorem~\ref{contimage} and the fact that $2^X$ is embedded into
$wL(\mathfrak{A})$. The set $\hat{\Sigma}^1_{-1}$ will then make
sure that $C'$ is in fact mapped onto $C$ by the map $f$.

Let us further construct the sets $\hat{\Sigma_n}$ for
$0<n<\omega$ in the same manner as we have constructed the set
$\Sigma_n$ in the previous section. So that if we have a model
$\mathfrak{A}$ of $\hat{\Sigma}$, the lattice $L(\mathfrak{A})$
will be a normal distributive and disjunctive lattice that models
the sentences~\ref{HI}, \ref{dim=1} and~$\conn(\1)$.

To prove the consistency of $\hat{\Sigma}$ it is enough to prove
the following lemma.
\begin{lemma}
For every finite $\Sigma'\in[\hat{\Sigma}]^{<\omega}$ there is a
metric continuum $X(\Sigma')$, and an interpretation function
$\mathcal{I}:\hat{K}\rightarrow 2^{X(\Sigma')}$ such that
$(2^{X(\Sigma')},\mathcal{I})$ is a model for $\Sigma'$.
\end{lemma}
\begin{proof}
Suppose we have a metric continuum $X(\Sigma'')$ for every subset
$\Sigma''$ of a given $\Sigma'\in[\Sigma]^{<\omega}$ such that
there is an interpretation function
$\mathcal{I}_{\Sigma''}:\hat{K}\rightarrow 2^{X(\Sigma'')}$ such
that $(2^{X(\Sigma'')},\mathcal{I}_{\Sigma''})\models\Sigma''$. We
want to show that there exists a metric continuum $X(\Sigma')$ and
an interpretation function
$\mathcal{I}_{\Sigma'}:\hat{K}\rightarrow 2^{X(\Sigma')}$ such
that $(2^{X(\Sigma')},\mathcal{I}_{\Sigma'})\models \Sigma'$.

Let $\theta$ be an $\sqsubset$-maximal sentence in $\Sigma'$ that
is of interest (see remark~\ref{wellorder}). If $\theta$ is an
element of $\hat{\Sigma}_0$, $\hat{\Sigma}_{5n+1}$,
$\hat{\Sigma}_{5n+2}$ or $\hat{\Sigma}_{5n+3}$ then we can choose
$X(\Sigma')=X(\Sigma'')$ and redefine the interpretation function
$\mathcal{I}$ in a natural way to obtain the wanted result. So let
us suppose that $\theta$ is an element of $\hat{\Sigma}_{-1}$,
$\hat{\Sigma}_{5n+4}$ or $\hat{\Sigma}_{5(n+1)}$ for some
$n<\omega$.

If $\theta$ is an element of $\hat{\Sigma}_{-1}$. Then, as
$\theta$ is the $\sqsubset$-maximal sentence in $\Sigma'$ of
interest, no $\phi\in\Sigma'$ is an element of
$\hat{\Sigma}_{5n+4}$ or $\hat{\Sigma}_{5(n+1)}$ for any
$n<\omega$. As $2^X$ is a normal, distributive and disjunctive
lattice and as $X$ is a continuum, we have that the lattice $2^X$
with an obvious interpretation function $\mathcal{I}$ is even a
model for $\triangle_{2^X}\cup\hat{\Sigma}_{-1}\cup\Sigma'$.

Suppose now that $\theta$ is an element of $\hat{\Sigma}_{5n+4}$
for some $n<\omega$. If we look at the construction in
subsection~\ref{theta4} we know that the function $(id\times h)$
is a closed monotone map from $X(\Sigma'')$ onto $X(\Sigma')$. So
inverse images of connected sets are connected and all the
sentences of $\hat{\Sigma}_{-1}$ in $\Sigma''$ that were true
(false) in the model $(2^{X(\Sigma'')},\mathcal{I}_{\Sigma''})$,
stay true (resp. false) in the model
$(2^{X(\Sigma')},\mathcal{I})$ as we get from
subsection~\ref{theta4}.

Finally, suppose that $\theta$ is an element of some
$\hat{\Sigma}_{5(n+1)}$. Lets take a look at the construction of
$X(\Sigma')$ in subsection~\ref{theta5}. Let $\pi$ be the map of
$X(\Sigma')$ onto $X(\Sigma'')$ as given in
subsection~\ref{theta5}. Consider the following lemma.
\begin{claim}
\label{preimconn}
    For every connected subset $A$ of $X(\Sigma'')$ there exists a
    connected set $C(A)\subset C=X(\Sigma')$ such that $\pi[C(A)]=A$.
\end{claim}
\begin{proof}
    Suppose we have $A\subset X(\Sigma'')$ connected. If we look at the
    image of $A$ under the function $f$ there are a number of
    possibilities:
    \begin{enumerate}
    \item
    \label{case1}
        $f[A]\subset [0,\frac{2}{3}]$ and
        $f[A]\cap[0,\frac{1}{3})\not=\emptyset$ ($f[A]\subset
        [\frac{1}{3},1]$ and
        $f[A]\cap(\frac{2}{3},1]\not=\emptyset$),
    \item
    \label{case2}
        $f[A]\subset[\frac{1}{3},\frac{2}{3}]$,
    \item
    \label{case3}
        $f[A]\setminus[0,\frac{2}{3}] \not=\emptyset \not=
        f[A]\setminus[\frac{1}{3},1]$
    \end{enumerate}
    In case~\ref{case1} we have $\{\frac{1}{4}\}\times A$
    ($\{\frac{3}{4}\}\times A$) is a connected subset of $X^+$
     which must intersect the component $C$, as every other
     component is mapped onto some subset of $X(\Sigma'')$, which is
     mapped into $[\frac{1}{3},\frac{2}{3}]$ by the function $f$.

     In case~\ref{case2} we have that the component $C$ must
     intersect at least one of the connected subsets
     $\{\frac{1}{4}\}\times A$, $\{\frac{1}{2}\}\times A$ or
     $\{\frac{3}{4}\}\times A$, as $C$ is mapped onto $X(\Sigma'')$,
     $X(\Sigma'')$ is connected and $f$ is continuous.

     In case~\ref{case3} we can, as above, assuming
     $A^+(=\pi^{-1}[A])=F+G$, construct closed and disjoint subsets
     $A^*$ and $B^*$ of $A$ which cover it. Again the image under
     $\pi$ is either all of $A$ or a proper subset of $A$. The
     (unique) piece
     that maps onto the whole of $A$ must intersect the set $C$,
     and so is contained in it.

     This ends the proof of the claim.
\end{proof}
We have that
$(2^{X(\Sigma')},\mathcal{I}_{\Sigma'})\models\Sigma'\setminus\hat{\Sigma}_{-1}$,
we now define a new interpretation function $\mathcal{I}$ on
$\hat{K}$ to $2^{X(\Sigma')}$ such that
$(2^{X(\Sigma')},\mathcal{I})$ will be a model for $\Sigma'$. Note
that the set of constants that are mentioned in the set $\Sigma'$
is a finite subset of $\hat{K}$, and let $\hat{K}(\Sigma')$ denote
this finite subset. We will define the interpretation under
$\mathcal{I}$ of the constants in $\hat{K}(\Sigma')$ 'from the
bottom up'.

By claim~\ref{preimconn} for every
$k_{-2,\alpha}\in\hat{K}(\Sigma')$ such that
$C=\mathcal{I}_{\Sigma''}(k_{-2,\alpha})$ is a connected subset of
$X(\Sigma'')$ we can find a connected subset $C'$ of $X(\Sigma')$
that maps onto $C$ by the map $\pi$. Let the $\mathcal{I}$
interpretation of the constant $k_{-2,\alpha}$ be this connected
set $C'$ in $X(\Sigma')$.

For all those $k$ in $\hat{K}(\Sigma')\cap\hat{K}_{-2}$ that have
no connected interpretation in $X(\Sigma')$ and  for all the
constants $k$ in $\hat{K}(\Sigma')\cap\hat{K}_{-1}$ the
interpretation under $\mathcal{I}$ will be the same as the
interpretation under $\mathcal{I}_{\Sigma'}$. So for those
$k\in\hat{K}(\Sigma')$ we have

$$\mathcal{I}(k_{-1,\alpha})=\mathcal{I}_{\Sigma'}(k_{-1,\alpha})=
\pi^{-1}[\mathcal{I}_{\Sigma''}(k_{-1,\alpha})]\cap X(\Sigma').$$

The interpretations of the rest of the constants in
$\hat{K}(\Sigma')$ will follow from the interpretations of the
constants we have just defined, because the interpretation of
every constant depends on just a finite set of other constants and
we just have to make sure that we define their interpretation in
the right order.

As $\mathcal{I}(k)\subset\mathcal{I}_{\Sigma'}(k)$ for all
$k\in\hat{K}(\Sigma')$ and for all
$k\in\hat{K}(\Sigma')\cap\hat{K}_{-2}$ such that
$\mathcal{I}_{\Sigma'}(k)$ is connected $\mathcal{I}(k)$ is also
connected we have that all the sentences of $\hat{\Sigma}_{-1}$
true (false) in $(2^{X(\Sigma'')},\mathcal{I}_{\Sigma''})$ are
true (false) in $(2^{X(\Sigma'')},\mathcal{I})$. The true or
falseness of the other sentences in $\Sigma'$ have not been
affected by the new interpretation function $\mathcal{I}$, and we
have completed the proof.
\end{proof}

%----------------------------------------------------------------
%
%   THE MA\'{C]kowiak-TymchATYN THEOREM
%
%----------------------------------------------------------------
\subsection{The Ma\'{c}kowiak-Tymchatyn theorem}
As we have seen in the previous section the set of sentences
$\hat{\Sigma}$ in the language $\{\linf,\lsup,\0,\1\}\cup\hat{K}$
is consistent. Let $\mathfrak{A}$ be a model for $\hat{\Sigma}$.
This model gives us a normal distributive and disjunctive lattice
$L(\mathfrak{A})$ which models the sentences~\ref{HI}, \ref{dim=1}
and~$\conn(\1)$. There also exists an embedding of $2^X$ into this
lattice $L(\mathfrak{A})$ (remember that we showed that, with
$K_{-2}$ an enumeration of the lattice $2^X$ the set $\Sigma$ is a
consistent set of sentences in the language
$\{\linf,\lsup,\0,\1\}\cup K$). All this implies that the Wallman
space $wL(\mathfrak{A})$, is a one-dimensional hereditarily
indecomposable continuum which admits a weakly confluent
surjection onto the metric continuum $X$.

Now we only have to make sure that there exists such a space that
is of countable weight to complete the proof of the Ma\'{c}kowiak
Tymchatyn theorem.
\begin{theorem}\cite{HvMP}
    Let $f:Y\rightarrow X$ be a continuous surjection between
    compact Hausdorff spaces. Then $f$ can be factored as $h\circ
    g$, where $Y\stackrel{g}{\rightarrow}
    Z\stackrel{h}{\rightarrow}X$ and~$Z$ has the same weight as
    $X$ and shares many properties with $Y$ (for instance, if $Y$
    is one-dimensional so is $X$ or if $Y$ is hereditarily
    indecomposable, so is $X$).
\end{theorem}
\begin{proof}
    Let $\mathcal{B}$ a minimal sized lattice-base for the closed
    sets of $X$, and identify it with its copy
    $\{f^{-1}[B]:B\in\mathcal{B}\}$ in $2^Y$. By the
    L\"{o}wenheim-Skolem theorem there is an elementary sublattice
    of $2^Y$, of the same cardinality as $\mathcal{B}$
    such that $\mathcal{B}\subset D\prec 2^Y$. The space
    $wD$ is as required.
\end{proof}
Applying this theorem to the space $wL(\mathfrak{A})$ and the
weakly confluent map $f:wL(\mathfrak{A})\rightarrow X$ we get a
one-dimensional hereditarily indecomposable continuum $wD$ which
admits a weakly confluent map onto the space $X$ and moreover the
weight of the space $wD$ equals the weight of the space $X$. This
is exactly what we were looking for.
%
%%%%%%%%%%%%%%%%%%%%%%%%%%%%%%%%%%%%%%%%%%%%%%%%%%%%%%%%%%%%
%                                                          %
%          A TOPOLOGICAL PROOF                             %
%                                                          %
%%%%%%%%%%%%%%%%%%%%%%%%%%%%%%%%%%%%%%%%%%%%%%%%%%%%%%%%%%%%
%
\section{A topological proof of the Ma\'{c}kowiak-Tymchatyn
theorem}\label{topproof}
%----------------------------------------------------------------
%
%   THE MA\'{C]kowiak-TymchATYN THEOREM (TOPOLOGICAL PROOF)
%
%----------------------------------------------------------------
\subsection{the Ma\'{c}kowiak-Tymchatyn theorem}
After the above proof was found we realized that it could be
transformed into a purely topological proof, which we shall now
describe.

Let $X$ be a metric continuum. We are going to define a inverse
sequence of metric continua with onto bonding maps $\{\langle
X_n,f_n\rangle:n<\omega\}$,
\begin{equation*}
X=X_0\stackrel{f_1}{\longleftarrow}
X_1\stackrel{f_2}{\longleftarrow}
\cdots\stackrel{f_n}{\longleftarrow}
X_n\stackrel{f_{n+1}}{\longleftarrow}\cdots,
\end{equation*}
in such a way that the inverse limit space $X_\omega$

$$X_\omega=\lim_\leftarrow\{\langle X_n,f_n\rangle:n<\omega\}$$

is a hereditarily indecomposable one-dimensional continuum of
weight $w(X)$ such that $\pi_0:X_\omega\rightarrow X$ is a weakly
confluent and onto. Here, for every $n<\omega$ the continuous
function $\pi_n$ is defined by $\pi_n=\operatorname{proj}_n\beperk
X_\omega:X_\omega\rightarrow X_n$, where
$\operatorname{proj}_n:\Pi_{n<\omega}X_n\rightarrow X_n$ is the
projection.

Let us furthermore define maps $f^n_m:X_n\rightarrow X_m$ for
$m<n$ as
\begin{equation*}
f^n_m=\left\{
    \begin{array}{l}
        f_{m+1}\circ f_{m+2}\circ\cdots\circ f_n,\ \text{if}\ m+1<n\\
        f_{m+1}\ \text{if}\ m+1=n.
    \end{array}
\right.
\end{equation*}
The following lemma is well known.
\begin{lemma}
The family of all sets of the form $\pi_n^{-1}(F)$, where $F$ is
an closed subset of the space $X_n$ and $n$ runs over a subset $N$
cofinal in $\omega$, is a base for the closed sets of the limit of
the inverse sequence $\{\langle X_n,f_n\rangle:n<\omega\}$.
Moreover, if for every $n<\omega$ a base $\mathcal{B}_n$ for the
closed sets of space $X_n$ is fixed, then the subfamily of those
$\pi_n^{-1}(F)$ for which $F\in\mathcal{B}_n$, also is a base for
the closed sets of $X_\omega$.
\end{lemma}

To make sure that the space $X_\omega$ is one-dimensional, it is
sufficient to show that $\{\pi_k^{-1}(F):F\in\mathcal{B}_k\
\text{and}\  k<\omega\}$ is a model of sentence~\ref{dim=1}.

Let $s:\omega\rightarrow\omega\times\omega$ be an onto map in such
a way that for every $n,m<\omega$ we have $s^{-1}(\langle
n,m\rangle)\geq\max\{n,m\}$. For instance, we take an onto map
$g:\omega\rightarrow\omega\times\omega\times\omega$ and given
$g(n)=\langle p,q,r\rangle$ we define $s(n)$ by
\begin{equation*}
s(n)=\left\{
        \begin{array}{l}
            \langle p,q\rangle\ \text{if}\ n\geq\max\{p,q\},\\
            \langle 0,0\rangle\ \text{otherwise}.
        \end{array}
\right.
\end{equation*}

Let $X_0=X$ and suppose we have defined the pairs $\langle
X_m,f_m\rangle$ and the bases $\mathcal{B}_m$ for every $m<n$. And
suppose that we have also defined an enumeration of all the
triples of $\mathcal{B}_m$ that have empty intersection for every
$m<n$. Let $\{G_k^m:k<\omega\}$ be an enumeration of the set
$\{G\in[\mathcal{B}_m]^3:\bigcap G=\emptyset\}$ for $m<n$, write
$G_k^m=\{a_k^m,b_k^m,c_k^m\}$.

The way we now define the space $X_n$ and the onto map
$f_n:X_n\rightarrow X_{n-1}$ will be as follows.

Suppose $s(n)=\langle k,m\rangle$, we consider the closed sets
$(f^{n-1}_k)^{-1}(a_m^k)$, $(f^{n-1}_k)^{-1}(b_m^k)$ and
$(f^{n-1}_k)^{-1}(c_m^k)$ of $X_{n-1}$. They have empty
intersection.

If there exist sets $x$, $y$ and~$z$ in $2^{X_{n-1}}$ such that
\begin{eqnarray*}
(f^{n-1}_k)^{-1}(a_m^k)\subset x, (f^{n-1}_k)^{-1}(b_m^k)\subset
y\ \text{and}\ (f^{n-1}_k)^{-1}(c_m^k)\subset z,\\ x\cap y\cap
z=\emptyset\ \text{and}\  x\cup y\cup z=X_{n-1},
\end{eqnarray*}
then we let $X_n=X_{n-1}$, $f_n=\id_{X_n}$ and we choose a
countable base $\mathcal{B}_n$ for the closed sets of $X_n$ such
that $\mathcal{B}_{n-1}\cup\{x,y,z\}\subset\mathcal{B}_n$.

If there do not exist such sets $x$, $y$ and~$z$ in $2^{X_n}$ then
we use the construction in subsection~\ref{theta4} to find a
(metric) continuum $X_n$ and a continuous onto map
$f_n:X_n\rightarrow X_{n-1}$, such that in $X_n$ there are closed
sets $x$, $y$ and~$z$ in $X_n$ such that
\begin{equation*}
    \begin{array}{l}
    (f^n_k)^{-1}(a^k_m)\subset x,\
    (f^n_k)^{-1}(b^k_m)\subset y,\
    (f^n_k)^{-1}(c^k_m)\subset z,\\
    x\cap y\cap z=\emptyset\ \text{and}\ x\cup y\cup z=X_n
    \end{array}
\end{equation*}
Let $\mathcal{B}_n$ be some countable base for the closed sets of
$X_n$ such that
$\{(f_n)^{-1}(F):F\in\mathcal{B}_{n-1}\}\subset\mathcal{B}_n$ and
$x,y,z\in\mathcal{B}_n$.

After we have chosen the base $\mathcal{B}_n$ we can choose some
enumeration of all the triples of $\mathcal{B}_n$ that have empty
intersection.

We do not get into trouble by considering base elements of some
base for the closed sets of $X_\omega$ which have not yet been
defined, because this will not happen by the way the function $s$
is defined and the bases $\mathcal{B}_n$ are chosen.

The limit $X_\omega(s)$ of the inverse sequence $\{\langle
X_n,f_n\rangle:n<\omega\}$ is a continuum, as all the spaces $X_n$
are continua, moreover, as the base
$\{\pi_n^{-1}(F^n_k):k,n<\omega\}$ of the space $X_\omega(s)$
models the sentence~\ref{dim=1} we have that $X_\omega(s)$ is
one-dimensional. As all the spaces $X_n$ are compact and all the
bonding maps $f_n$ are onto, we have that
$\pi_0:X_\omega(s)\rightarrow X$ is a continuous onto map.

In a similar way we can construct a function
$t:\omega\rightarrow\omega\times\omega\times\omega$, an onto map
in such a way that for all $k,l,m<\omega$ we have $t^{-1}(\langle
k,l,m\rangle)\geq\max\{k,l,m\}$, and use it together with the
construction in subsection~\ref{theta5} to define, given $X_0=X$,
$X_n$ and $f_n$ so that $X_\omega(t)$, the inverse limit of the
sequence $\{\langle X_n,f_n\rangle:n<\omega\}$ is a hereditarily
indecomposable continuum which admits a continuous onto map,
$\pi_0$ onto the space $X$.

We can combine these two constructions by defining the function
$r$ by letting $r(2n)$ equal $s(n)$ and $r(2n+1)$ equal $t(n)$ for
every $n<\omega$. Define $X_0=X$ and use the construction in
subsection~\ref{theta4} if $n$ is even and the construction in
subsection~\ref{theta5} if $n$ is odd to construct $X_n$ and
$f_n$.

Let $X_\omega(r)$ be the inverse limit of the inverse sequence
$\{\langle X_n,f_n\rangle\}$ we have constructed with the aid of
the function $r$ as described in the previous paragraph.

As $\mathcal{B}=\{\pi_n^{-1}(F^n_k):n,k<\omega\}$ is a base for
the closed sets of $X_\omega(r)$ we see that
$w(X_\omega(r))=w(X)=\aleph_0$. The space $X_\omega(r)$ is a
one-dimensional hereditarily indecomposable continuum as, by
construction $\mathcal{B}$ is a model of the sentences~\ref{dim=1}
and~\ref{HI}. So by the following claim we have proven the
Ma\'{c}kowiak-Tymchatyn theorem.
\begin{claim}
The map $\pi_0$ is a weakly confluent map from $X_\omega(r)$ onto
$X$.
\end{claim}
\begin{proof}
Suppose we have a subcontinuum $C$ of the space $X$, we want to
find a subcontinuum $C'$ of $X_\omega(r)$ such that $\pi_0[C']=C$.
As, by construction, the $f_{2n}$'s are monotone closed maps from
$X_{2n}$ onto $X_{2n-1}$ and the $f_{2n+1}$'s are weakly
confluent, we can define an inverse sequence $\{\langle Y_n,
g_n\rangle:n<\omega\}$ such that $Y_0=C$ and $Y_n$ is, for every
$n$, some subcontinuum of $f_n^{-1}(Y_{n-1})$ that is mapped onto
$Y_{n-1}$ by the map $f_n$ and the map $g_n$ is the restriction of
the map $f_n$ to the subspace $Y_n$ of $X_n$. Let $C'$ be the
inverse limit of the inverse sequence $\{\langle
Y_n,g_n\rangle:n<\omega\}$, so

$$C'=\lim_\leftarrow\{\langle Y_n,g_n\rangle:n<\omega\}.$$

We have that $C'$ is a closed subspace of the space $X_\omega(r)$.
Furthermore it is a continuum as it is an inverse limit of
continua, so it is a subcontinuum of $X_\omega(r)$. As $\pi_n$
maps $C'$ onto the $Y_n$'s we have proven the claim.
\end{proof}
%
%----------------------------------------------------------------
%
%   THE EXTENDED MA\'{C]kowiak-TymchATYN THEOREM
%
%----------------------------------------------------------------
\subsection{The extended Ma\'{c}kowiak-Tymchatyn theorem}
Given a continuum $X$ we will construct an inverse sequence
$\{\langle X_\alpha,f_\alpha\rangle:\alpha<w(X)\}$ such that the
inverse limit space $Y$ is a hereditarily indecomposable continuum
of weight $w(X)$ and $\operatorname{dim}(Y)=1$ and there exists a
weakly confluent map of the space $Y$ onto~$X$. This is a somewhat
different proof than is given in the paper of Hart, Van Mill and
Pol (see~\cite{HvMP}).
\begin{equation*}
X=X_0\stackrel{f_1}{\longleftarrow}
X_1\stackrel{f_2}{\longleftarrow}
\cdots\stackrel{f_\alpha}{\longleftarrow}
X_{\alpha}\stackrel{f_{\alpha+1}}{\longleftarrow}\cdots\quad(\alpha<w(X)).
\end{equation*}
We are going to make sure that every $X_\alpha$ is a continuum of
weight $w(X)$ and that there exist some base $\mathcal{B}_\alpha$
for the closed sets of $X_\alpha$ of cardinality $w(X)$ such that
the base $\{\pi^{-1}(B):B\in\mathcal{B}_\alpha,\alpha<w(X)\}$ for
the closed sets of the space $Y$ will show that $Y$ is the desired
continuum.

For $\beta<w(X)$ a limit ordinal we let $X_\beta$ be the inverse
limit of the sequence $\{\langle
X_\gamma,f_\gamma\rangle:\gamma<\beta\}$, and we let
$\mathcal{B}_\beta$ be the set
$\{(\pi_\alpha^\beta)^{-1}(B):B\in\mathcal{B}_\alpha,\
\alpha<\beta\}$. This is a base for the closed sets of $X_\beta$
and $\vert\mathcal{B}_\beta\vert\leq w(X)$. Furthermore $X_\beta$
is a continuum as it is an inverse limit of continua.

Suppose we have defined the continua $X_\beta$ for
$\beta\leq\alpha$ for some $\alpha<w(X)$, as well as the bases
$\mathcal{B}_\beta$ for the closed sets of these spaces and for
every $\beta\leq\alpha$ we also have defined an enumeration
$\{G_\tau^\beta:\tau<\Gamma_\beta\}$ of the triples of elements of
$\mathcal{B}_\beta\setminus\{\emptyset\}$ that have empty
intersection, we write
$G_\tau^\beta=\{a_\tau^\beta,b_\tau^\beta,c_\tau^\beta\}$. Here
$\Gamma_\beta$ is some ordinal number less than or equal to
$w(X)$.

As in the previous section we can find a function
$s:w(X)\rightarrow w(X)\times w(X)$ such that for every
$\alpha,\beta< w(X)$ we have $s^{-1}(\langle
\alpha,\beta\rangle)\geq\max\{\alpha,\beta\}$. To find
$X_{\alpha+1}$ and $f_\alpha$ we do almost the same thing as we
have done in the previous section. If $s(\alpha)=\langle
\beta,\gamma\rangle$ we consider the closed sets
$a=f_\beta^\alpha(a_\gamma^\beta)$,
$b=f_\beta^\alpha(b_\gamma^\beta)$ and
$c=f_\beta^\alpha(c_\gamma^\beta)$ of the space $X_\alpha$.

If there exist $x$, $y$ and~$z$ in $2^{X_\alpha}$ such that
$a\subset x$, $b\subset y$, $c\subset z$, $x\cap y\cap
z=\emptyset$ and $x\cup y\cup z=X_\alpha$ then we let
$X_{\alpha+1}=X_\alpha$ and $f_{\alpha+1}=\id_{X_{\alpha+1}}$.

If there are no such $x$,$y$ and~$z$ in $2^{X_\alpha}$ then we
will do as in subsection~\ref{theta4}, but as in that section we
used a metric for $X$ we have to slightly alter the proof there.
As $X_\alpha$ is normal and $a\cap b\cap c=\emptyset$ we can find
a continuous function $f_a:X_\alpha\rightarrow [0,1]$ such that
$f_a(a)\subset\{0\}$ and $f_a(b\cap c)\subset\{1\}$. Now, as
$f_a^{-1}(\{0\})\cap b\cap c=\emptyset$ we can find a continuous
function $f_b:X_\alpha\rightarrow [0,1]$ such that
$f_b(b)\subset\{0\}$ and $f_b(f_a^{-1}(\{0\})\cap c)\subset\{1\}$.
Finally, since $f_a^{-1}(\{0\})\cap f_b^{-1}(\{0\})\cap
c=\emptyset$ we can find a continuous function
$f_c:X_\alpha\rightarrow [0,1]$ such that $f_c(c)\subset\{0\}$ and
$f_c(f_a^{-1}(\{0\})\cap f_b^{-1}(\{0\}))\subset\{1\}$. Now define
the function $f:X_\alpha\rightarrow\mathbb{R}^3$ by

$$f(x)=(\kappa_a(x),\kappa_b(x),\kappa_c(x)),$$

where $\kappa_a:X_\alpha\rightarrow [0,1]$ is defined by

$$\kappa_a(x)=\frac{f_a(x)}{f_a(x)+f_b(x)+f_c(x)},$$

and $\kappa_b$ and $\kappa_c$ are likewise defined. The function
$f$ maps $X_\alpha$ into the triangle that is the convex hull of
the points $\{(0,0,1),(0,1,0),(1,0,0)\}$ in $\mathbb(R)^3$ just as
in subsection~\ref{theta4} and from this point on we can follow
the method in subsection~\ref{theta4} to find a continuum
$X_{\alpha+1}$ and a continuous onto map
$f_{\alpha+1}:X_{\alpha+1}\rightarrow X_\alpha$ such that there
exist $x,y,z\in 2^{X_{\alpha+1}}$ such that
$f_{\alpha+1}^{-1}(a)\subset x$, $f_{\alpha+1}^{-1}(b)\subset y$
and $f_{\alpha+1}^{-1}(c)\subset z$, $x\cap y\cap z=\emptyset$ and
$x\cup y\cup z= X_{\alpha+1}$. Now let $\mathcal{B}_{\alpha+1}$ be
a base for the closed sets of $X_{\alpha+1}$ such that
$\{(f_{\alpha+1})^{-1}(B):B\in\mathcal{B}_\alpha\}\cup
\{x,y,z\}\subset\mathcal{B}_{\alpha+1}$ and
$\vert\mathcal{B}_{\alpha+1}\vert\leq w(X)$. Enumerate the set of
triples of $\mathcal{B}_{\alpha+1}\setminus\{\emptyset\}$ with
empty intersection as
$\{G_\tau^{\alpha+1}:\tau<\Gamma_{\alpha+1}\}$, where
$\Gamma_{\alpha+1}$ is some ordinal number less than or equal to
$w(X)$.

In a similar way we can find an (transfinite) inverse sequence
such that the inverse limit is a hereditarily indecomposable
continuum of the same weight as $X$ and for which the map $\pi_0$
is a continuous onto map between the limit and the space $X$.

As in the previous section we can combine these two (we take care
of the hereditary indecomposability at even ordinal stages and we
take care that the dimension of the limit space will not exceed
one at the odd ordinal stages), to find a transfinite inverse
sequence such that the inverse limit is a one-dimensional
hereditarily indecomposable continuum that admits a continuous map
onto the space $X$.  After some thought, as in the previous
section we see that this continuous map is in fact a weakly
confluent map.
%
%%%%%%%%%%%%%%%%%%%%%%%%%%%%%%%%%%%%%%%%%%%%%%%%%%%%%%%%%%%%%
%                                                           %
%                         THE BIBLIOGRAPHY                  %
%                                                           %
%%%%%%%%%%%%%%%%%%%%%%%%%%%%%%%%%%%%%%%%%%%%%%%%%%%%%%%%%%%%%

\end{document}